\documentclass[12pt]{amsart}
\pdfoutput=1

\usepackage[utf8]{inputenc}
\usepackage{lmodern}
\usepackage[T1]{fontenc}
\usepackage[leqno]{amsmath}
\allowdisplaybreaks
\usepackage{amssymb} % also loads amsfonts
\usepackage{amsthm}
\usepackage{mathtools}
\usepackage[usenames,dvipsnames,svgnames,table]{xcolor}
\usepackage{embrac} % upright parentheses inside \emph{}
\ChangeEmph{(}[-.01em,.04em]{)}[.04em,-.05em]
\usepackage{microtype}
\usepackage{textcase}
\usepackage{enumitem}
\usepackage{graphicx}
\usepackage{float}
\usepackage{booktabs}
\usepackage{csquotes}
\usepackage{multirow}
\usepackage{tikz}
\usepackage{subcaption}
\tikzset{node distance=2cm, auto}
\usetikzlibrary{angles,quotes,calc,positioning,shapes,matrix,arrows}
\usepackage{tikz-cd}
\usepackage{tkz-euclide}
\usepackage{etoolbox}
\usepackage{url}
\usepackage[noadjust]{cite}
\usepackage[hyperindex=true, colorlinks=true, linkcolor=Maroon, citecolor=DarkGreen, urlcolor=NavyBlue, breaklinks=true,linktocpage=true,linktoc=all]{hyperref}
\usepackage[english,capitalise,noabbrev]{cleveref}
\crefdefaultlabelformat{#2\textup{#1}#3}
\crefname{equation}{equation}{equations}
\newcommand{\N}{\mathbb{N}}
\newcommand{\Z}{\mathbb{Z}}

\newcommand{\R}{\mathbb{R}}
\newcommand{\C}{\mathbb{C}}

\renewcommand{\S}{\mathbb{S}}

\newcommand{\RP}{\mathbb{RP}}

\newcommand{\dd}{\mathop{}\!\mathrm{d}}

\newcommand{\contained}{\subset}

\newcommand{\union}{\cup}
\newcommand{\bigunion}{\bigcup}

\newcommand{\suchthat}{\, : \,}
\newcommand{\st}{\suchthat}

\newcommand{\defined}{\coloneqq}

\newcommand{\bigdirectsum}{\bigoplus}

\renewcommand{\O}{\operatorname{O}}

\newcommand{\SO}{\operatorname{SO}}

\newcommand{\Id}{I}

\newcommand{\rank}{\operatorname{rank}}

\newcommand{\from}{\colon}
 
\newcommand{\composed}{\circ}

\DeclarePairedDelimiterX\norm[1]\lVert\rVert{\ifblank{#1}{\:\cdot\:}{#1}}
\DeclarePairedDelimiterX\abs[1]\lvert\rvert{\ifblank{#1}{\:\cdot\:}{#1}}
\DeclarePairedDelimiterX\set[1]{\{}{\}}{\ifblank{#1}{\: \:}{#1}}
\DeclarePairedDelimiterX\innerprod[2]\langle\rangle
{\ifblank{#1#2}{\,\cdot,\cdot\,}
	{#1,#2}}
\DeclarePairedDelimiterX\floor[1]\lfloor\rfloor{\ifblank{#1}{\:\cdot\:}{#1}}
\DeclarePairedDelimiterXPP\expectation[2]{\ifblank{#1}{\mathbb{E}}{\mathbb{E}_{#1}}}\lbrack\rbrack{}{\ifblank{#2}{\:\cdot\:}{#2}}
\newcommand{\Grass}[3][]{\ifblank{#1}{\operatorname{Gr}_{#2,#3}}
	{\operatorname{Gr}_{#2,#3}(#1)}}

\newcommand{\transpose}{^{\text{\tiny$\mathsf{T}$}}}
\newcommand{\length}{\ell}
\newcommand{\ParamGrass}{\mathcal{P}}
\newcommand{\Pol}{\operatorname{Pol}}
\newcommand{\PolFixed}{\Pi} % space of spherical harmonics of exact degree
\newcommand{\DegreePart}[1]{\abs{#1}}
\newcommand{\restricted}[1][\big]{#1|}
\newcommand{\Span}{\operatorname{span}}

\theoremstyle{plain}
\newtheorem{theorem}{Theorem}[section]

\newtheorem{lemma}[theorem]{Lemma}
\newtheorem{proposition}[theorem]{Proposition}

\theoremstyle{definition}

\newtheorem{definition}[theorem]{Definition}

\theoremstyle{remark}

\begin{document}

	\title[Asymptotically optimal $t$-design curves on $\operatorname{Gr}_{2,4}$ and $\operatorname{SO}(3)$]{Asymptotically optimal $t$-design curves\\ on $\operatorname{Gr}_{2,4}$ and $\operatorname{SO}(3)$}
	
    \author{Martin Ehler}
	\address{Martin Ehler: Fakultät für Mathematik, Universität Wien, Oskar-Morgenstern-Platz 1, 1090 Vienna, Austria}
    \email{martin.ehler@univie.ac.at}

    \author{Clemens Karner}
	\address{Clemens Karner: Fakultät für Mathematik, Universität Wien, Oskar-Morgenstern-Platz 1, 1090 Vienna, Austria}
	\email{clemens.karner@univie.ac.at}
    
	\author{Pedro R. López-Gómez}
	\address{Pedro R. López-Gómez: Departamento de Matemáticas, CUNEF Universidad,  C/ Almansa, 101, 28040 Madrid, Spain}
	\email{pedro.lopezgomez@cunef.edu}
	
	\date{\today{}}
	
	\thanks{Clemens Karner has been supported by funding of the Austrian Science Fund, project T1307-N ``Image based Data Evaluation Analyses and Medical Application''. Pedro R. López-Gómez has been supported by grant PID2020-113887GB-I00 funded by MCIN/AEI/10.13039/501100011033 and by grant PRE2021-097772 funded by MCIN/AEI/10.13039/501100011033 and by ``ESF Investing in your future''.}
	
	

	\begin{abstract}
		In analogy with spherical $t$-design points, we study $t$-design curves on the Grassmannian $\operatorname{Gr}_{2,4}$ and on the rotation group $\operatorname{SO}(3)$, that is, closed curves along which the normalized line integral of every polynomial of degree at most $t$ equals its average over the manifold. Using the identifications of $\operatorname{Gr}_{2,4}$ with $(\mathbb{S}^2\times\mathbb{S}^2)/\{\pm 1\}$ and of $\operatorname{SO}(3)$ with $\mathbb{S}^3/\{\pm 1\}$, we construct sequences of $t$-design curves whose lengths are of order $t^3$ and $t^2$, respectively, which is asymptotically optimal.
	\end{abstract}
	
	\maketitle

%============================================================

\section{Introduction}

A \emph{$t$-design point configuration} on the $d$-sphere $\S^d=\{x\in \R^{d+1}\st \norm{x}=1\}$ is a finite set of points $X=\{x_1,\dotsc, x_n\}$, with $n \in \N$, such that
\begin{equation*}
    \frac{1}{n}\sum_{i=1}^{n}f(x_i)=\int_{\S^d} f\quad\text{for all $f\in\Pol_t(\S^d)$},
\end{equation*}
where $\Pol_t(\S^d)$ is the set of polynomials in $d+1$ variables of degree at most $t$ restricted to $\S^d$, the integral is taken with respect to the normalized surface measure, and $d,t\in \N$.  

Spherical $t$-design points have been studied extensively; see, for instance, \cite{tdesignpoints-hardin1996mclaren,tdesignpoints-reznick1995some,tdesignpoints-seidel2001definitions,tdesignpoints-sloane2003spherical,tdesignpoints-womersley2018efficient,tdesignpoints-bondarenko2013optimal,tdesignpoints-boyvalenkov2001uniqueness,tdesignpoints-goethals1981football,tdesignpoints-graf2011computation}. Further existence results on more general compact Riemannian manifolds can be found in \cite{tdesignpoints-gariboldi2021optimal,tdesignpoints-etayo2018asymptotically,tdesignpoints-ehler2021asymptotically}.

The concept of spherical designs was introduced in 1977 by Delsarte, Goethals, and Seidel in their fundamental paper \cite{tdesignpoints-delsarte1977spherical}. Denoting by $N(d,t)$ the minimal number of points of a spherical $t$-design on $\S^d$ for each $d,t\in\N$, they proved the lower bound $N(d,t)\geq c_d t^d$. Here and throughout this text, $c_d$ and $C_d$ denote sufficiently small and sufficiently large positive constants, respectively, depending only on $d$. Several years later, in 1984, Seymour and Zaslavsky \cite{tdesignpoints-seymour1984averaging} proved the existence of $t$-design points on $\S^d$ for all $d,t\in \N$. A key open question was the existence of a tight upper bound on the number of points. The first results in this regard were obtained in 1991 by Wagner and Volkmann \cite{tdesignpoints-wagner1991averaging}, who showed that $N(d,t) \leq \cramped{C_d t^{C_d d^4}}$, and in 1992 by Bajnok \cite{tdesignpoints-bajnok1992construction}, who obtained $N(d,t) \leq \cramped{C_d t^{C_d d^3}}$. Shortly afterwards, in 1993,  Korevaar and Meyers proved in their well-known paper \cite{tdesignpoints-korevaar1993spherical} that $N(d,t)\leq \cramped{C_dt^{(d^2+d)/2}}$ and conjectured that $N(d,t)\leq C_dt^d$, which became known as the Korevaar--Meyers conjecture. Sequences of spherical $t$-designs $(X_t)_{t\in\N}$ with $\abs{X_t}\asymp t^d$ are called \emph{asymptotically optimal}\footnote{The notation $f(t)\asymp g(t)$ means that $f(t)=O(g(t))$ and $g(t)=O(f(t))$, that is, $f$ and $g$ are asymptotically of the same order.}. In this regard, one of the major breakthroughs was the proof by Bondarenko, Radchenko, and Viazovska \cite{tdesignpoints-bondarenko2013optimal} in 2013 of the aforementioned Korevaar--Meyers conjecture on the existence of asymptotically optimal $t$-design configurations on $\S^d$, a problem that had been open for 20 years. Their result asserts that, for every $d$, there exists a sequence of spherical $t$-designs $(X_t)_{t\in\N}$ on $\S^d$ with $\abs{X_t} \asymp t^d$.

\subsection*{Spherical $t$-design curves}

In 2023, the concept of $t$-design points was generalized to curves on the $d$-sphere in \cite{tdesigncurves-EhlerGrochenig2023}. A \emph{$t$-design curve} on $\S^d$ is a closed, piecewise smooth curve $\gamma\from [0,1] \to \S^d$ of length $\length(\gamma)=\int_0^1 \norm{\dot{\gamma}(s)} \dd s$ such that the line integral $\int_{\gamma}f=\int_{0}^{1}f(\gamma(s))\norm{\dot{\gamma}(s)} \dd s$ integrates all polynomials in $\Pol_t(\S^d)$ exactly, that is,
\begin{equation*}
		\frac{1}{\length(\gamma)}\int_{\gamma}f=\int_{\S^d}f\quad\text{for all $f\in\Pol_t(\S^d)$}.
\end{equation*}
In contrast to \cite{tdesigncurves-EhlerGrochenig2023}, we allow arbitrary self-intersections. It is shown in \cite{tdesigncurves-EhlerGrochenig2023} that the length of every $t$-design curve is bounded from below by $c_dt^{d-1}$, with a constant $c_d>0$ that may depend on the dimension $d$ but is independent of $t$ and $\gamma$. This lower bound was recently made explicit by Ehler and Gröchenig \cite{ehlergroechenig2026quantitative} and sharpened by King and Mixon \cite{king2026t2short}; see \cref{sec:sphere}.

Explicit examples of spherical $t$-design curves can be found in \cite{tdesigncurves-EhlerGrochenig2023,tdesigncurves-ehler2025geodesic,king2026t2short}. An analogue of the Korevaar--Meyers conjecture for the case of $t$-design curves asks for the existence of asymptotically optimal $t$-design curves on $\S^d$, that is, a sequence $(\gamma_t)_{t\in\N}$ of $t$-design curves whose length grows at most like $t^{d-1}$. This was settled for the $2$-sphere in \cite{tdesigncurves-EhlerGrochenig2023} and for the $3$-sphere in \cite{tdesigncurves-lindblad2024asymptotically}, with sequences of $t$-design curves of lengths $\length(\gamma_t)\asymp t$ and $\length(\gamma_t)\asymp t^2$, respectively. Generalizations of this question to approximate and weighted $t$-design curves on $\S^d$ have recently been studied in \cite{tdesigncurves-lindblad2025generalizations}. Beyond spheres, Lindblad constructed in \cite{tdesigncurves-lindblad2024asymptotically} asymptotically optimal $t$-design curves on the torus $\mathbb{T}^d=\R^d/\Z^d$ for trigonometric polynomials of total degree at most $t$, and Kunis \cite{kunis2026Td} derived lower bounds of the same order on the length of such curves.

\subsection*{$t$-design curves on $\Grass{2}{4}$ and $\SO(3)$}

We are especially interested in asymptotically optimal $t$-design curves on the Grassmann manifold $\Grass{2}{4}$ and on the rotation group $\SO(3)$. Both manifolds carry natural analogues of the polynomial spaces $\Pol_t(\S^d)$: on $\Grass{2}{4}$, viewed as the set of rank-$2$ orthogonal projectors on $\R^4$, these are the restrictions to the manifold of polynomials of degree at most $t$ in the matrix entries, and on $\SO(3)$ they are the span of the Wigner $D$-functions $D^n_{k,k'}$ with $n\leq t$.

For $t$-design points, asymptotically optimal configurations on both manifolds exist by the work of Etayo, Marzo, and Ortega-Cerdà \cite{tdesignpoints-etayo2018asymptotically}, who proved the Korevaar--Meyers conjecture on compact algebraic manifolds, a class that includes both the Grassmannians and $\SO(3)$. To the best of our knowledge, no corresponding result for $t$-design curves on the Grassmannian was known before the present work.

The Grassmannian $\Grass{2}{4}$ encodes orthogonal projections for representing data in lower dimensions, with applications including the visualization of high-dimensional data sets \cite{applications-tSNE,applications-breger2024visclust} and the classification of clusters \cite{applications-breger2024visclust}. The rotation group $\SO(3)$ contains the rotations in three-dimensional space, with applications in texture analysis, where distributions of crystal orientations are modeled on $\SO(3)$ \cite{applications-so3}. Furthermore, \cite{hielscher2013kernel} treats kernel density estimation and \cite{hielscher2013funk} the numerical inversion of the Funk transform, both on $\SO(3)$. Fast Fourier transforms on $\SO(3)$ and on $\Grass{2}{4}$ were derived in \cite{DickEhlerGraefKrattenthaler2023}, based on the spectral decomposition of discrepancy kernels on these two manifolds. In \cite{EhlerGraefNeumayerSteidl2021}, measures on $\SO(3)$ and on $\Grass{2}{4}$ were approximated by curve-supported measures, by numerically minimizing the discrepancy between a given measure and the pushforward of the Lebesgue measure on $[0,1]$ under a Lipschitz curve. This is closely related to the notion of $t$-design curves investigated here.

\subsection*{Results}
In this paper, we construct asymptotically optimal $t$-design curves on these two manifolds. As the lower bound in \cref{thm:curve-lowerbound} carries over to compact Riemannian manifolds, and as $\Grass{2}{4}$ and $\SO(3)$ have dimensions $4$ and $3$, every sequence of $t$-design curves on these manifolds satisfies $\length(\gamma_t)\geq ct^3$ and $\length(\gamma_t)\geq ct^2$, respectively, for some constant $c>0$ \cite{tdesigncurves-EhlerGrochenig2023}. Here, asymptotically optimal means that the length of the curves matches this expected optimal length $t^{d-1}$ on a $d$-dimensional manifold, in accordance with the curve analogue of the Korevaar--Meyers conjecture. Our main results are the following:
\begin{enumerate}[label={(\Alph*)}]
	\item There exists a sequence $(\gamma_t)_{t\in\N}$ of $t$-design curves on $\Grass{2}{4}$ with $\length(\gamma_t)\asymp t^3$ (\cref{thm:optimal-curves-G24}).
	\item The construction of a sequence $(\gamma_t)_{t\in\N}$ of $t$-design curves on $\SO(3)$ with $\length(\gamma_t)\asymp t^2$ (\cref{thm:optimal-curves-SO3}).
\end{enumerate}
 Thus, the sequences constructed in (A) and (B) match these lower bounds and are therefore asymptotically optimal. For $\RP^3\cong\SO(3)$, this was mentioned by Lindblad in \cite[Section 3]{tdesigncurves-lindblad2024asymptotically}; \cref{thm:optimal-curves-SO3} gives a construction in terms of Wigner $D$-functions. A preliminary version of (A) appeared in \cite{sampta-gr24}.

\subsection*{Outline of the paper}

\Cref{sec:sphere} collects known results on $t$-design curves on the sphere. \Cref{sec:grassmann} introduces polynomials and $t$-designs on the Grassmannian together with the parametrization of $\Grass{2}{4}$ as $(\S^2\times\S^2)/\{\pm 1\}$. In \cref{sec:grassmann-tdesigns}, we define $t$-design curves on the Grassmannian and construct asymptotically optimal ones on $\Grass{2}{4}$. An explicit example of a $2$-design curve is presented in \cref{sec:grassmann-example}. In \cref{sec:so3-harmonic}, we develop the harmonic decomposition of $\SO(3)$ and its identification with $\S^3/\{\pm 1\}$, which \cref{sec:SO3} exploits to construct asymptotically optimal $t$-design curves on the rotation group. Finally, \cref{sec:other} extends the construction to the projective spaces $\RP^{2}\cong\Grass{1}{3}$ and $\RP^{3}\cong\Grass{1}{4}$ (and hence to $\Grass{2}{3}$ and $\Grass{3}{4}$) and recalls asymptotically optimal $t$-design curves on the torus $\mathbb{T}^d$.

\section{$t$-design curves on the sphere}\label{sec:sphere}

In this section, we collect known results on $t$-design curves on the sphere that we build upon in this paper. When passing from $t$-design points to curves on the $d$-sphere, the normalized counting measure is replaced by the normalized arclength measure along the curve, and the number of points by the length of the curve. Since a curve is a one-dimensional subset of the $d$-dimensional sphere, the analogue of the lower bound $N(d,t)\geq c_dt^d$ for spherical $t$-designs one expects is a lower bound of order $t^{d-1}$ on the length of $t$-design curves. Indeed, in \cite{tdesigncurves-EhlerGrochenig2023} it is shown that the length of every $t$-design curve is bounded from below as follows.

\begin{theorem}[{\cite[Theorem 1.1]{tdesigncurves-EhlerGrochenig2023}}]\label{thm:curve-lowerbound}
	Assume that a closed, piecewise smooth curve $\gamma\from [0,1]\to\S^d$ satisfies 
	\begin{equation*}
		\frac{1}{\length(\gamma)}\int_{\gamma}f=\int_{\S^d}f\quad\text{for all $f\in\Pol_t(\S^d)$}.
	\end{equation*}
	Then, its length is bounded from below by
	\begin{equation*}
		\length(\gamma)\geq c_dt^{d-1},
	\end{equation*}
	with some constant $c_d>0$ that may depend on the dimension $d$ but is independent of $t$ and $\gamma$.
\end{theorem}

Recently, this lower bound has been made explicit in \cite{ehlergroechenig2026quantitative}:
\begin{theorem}[{\cite[Theorem~5.2 and Corollary~5.3]{ehlergroechenig2026quantitative}}]\label{thm:ehlergroechenig-lowerbound}
	Let $t\geq d\geq 2$. Every $2t$-design curve $\gamma_{2t}$ on $\S^d$ satisfies
	\begin{equation*}
		\length(\gamma_{2t})\geq\frac{\sqrt{d+2}}{4}\,\frac{\dim\Pol_t(\S^d)}{t}\geq \frac{\sqrt{d+2}}{2d!}t^{d-1}.
	\end{equation*}
\end{theorem}
In particular, this provides an explicit value of the constant $c_d$ in \cref{thm:curve-lowerbound}. King and Mixon obtained a different explicit bound, valid for all $t$, by spectral methods:

\begin{theorem}[{\cite[Theorem 4]{king2026t2short}}]\label{thm:explicit-lowerbound}
	Every $t$-design curve $\gamma_t$ on $\S^d$ satisfies $\length(\gamma_t)\geq 2\pi\sqrt{A_{t,d}/B_{t,d}}$, where $A_{t,d}$ and $B_{t,d}$ are explicit constants depending only on $t$ and $d$. In particular,
	\begin{equation*}
		\length(\gamma_t)\geq \frac{\pi}{2^{d-2}d!}\sqrt{\frac{d+2}{3}}\, t^{d-1}.
	\end{equation*}
\end{theorem}

The bound $2\pi\sqrt{A_{t,d}/B_{t,d}}$ is attained for $t=1$ in all spheres and for $t=2$ in every odd-dimensional sphere, which, as noted in \cite{king2026t2short}, yields the first exact optimality results for $t$-design curves with $t>1$. The constant in front of $t^{d-1}$ in \cref{thm:explicit-lowerbound} improves on the one obtained from \cref{thm:ehlergroechenig-lowerbound} in every dimension.

In view of \cref{thm:curve-lowerbound}, a sequence of $t$-design curves on $\S^d$ is asymptotically optimal if its length grows at most like $t^{d-1}$. Such sequences of curves are known to exist on the $2$-sphere, as presented in \cite{tdesigncurves-EhlerGrochenig2023}, and on the $3$-sphere, by Lindblad \cite{tdesigncurves-lindblad2024asymptotically}:

\begin{theorem}[{\cite[Theorem 1.2]{tdesigncurves-EhlerGrochenig2023}}]\label{thm:optimal-curves-S2}
	In $\S^2$, there exists a sequence of $t$-design curves $(\gamma_t)_{t\in\N}$ with length $\length(\gamma_t)\asymp t$.
\end{theorem}

\begin{theorem}[{\cite[Theorem 1.2]{tdesigncurves-lindblad2024asymptotically}}]\label{thm:optimal-curves-S3}
	In $\S^3$, there exists a sequence of $t$-design curves $(\gamma_t)_{t\in\N}$ with length $\length(\gamma_t)\asymp t^2$.
\end{theorem}

We build on both results in the following sections. In \cref{sec:grassmann-tdesigns}, for the Grassmannian $\Grass{2}{4}$, we use \cref{thm:optimal-curves-S2} for the $2$-sphere in conjunction with an analogous theorem for $t$-design points, and, in \cref{sec:SO3}, for the rotation group $\SO(3)$, we employ \cref{thm:optimal-curves-S3} for the $3$-sphere.

\section{Harmonic decomposition of the Grassmannian}\label{sec:grassmann}

The (real, unoriented) Grassmann manifold $\Grass{m}{d}$ is defined as the set of $m$-dimensional subspaces of $\R^d$. It is a compact, connected Riemannian manifold of dimension $m(d-m)$ and a homogeneous space under the action of the orthogonal group $\O(d)$:
\begin{equation*}
    \Grass{m}{d}=\O(d)/\bigl(\O(m)\times \O(d-m)\bigr).
\end{equation*}
We will denote the unique $\O(d)$-invariant probability measure on $\Grass{m}{d}$ induced by the Haar measure on $\O(d)$ as $\sigma_{\Grass{m}{d}}$ and we will refer to it as the \emph{uniform measure} on $\Grass{m}{d}$.

In this work, we will use the following identification of $\Grass{m}{d}$ with the set of orthogonal projectors on $\R^d$ of rank~$m$:
\begin{equation*}
    \Grass{m}{d}=\set{P\in\R^{d\times d}\st P\transpose=P,\ P^2=P,\ \rank(P)=m}.
\end{equation*}

In this section, we recall the harmonic analysis on the Grassmannian $\Grass{m}{d}$ that underlies the definition of polynomials on this space, as well as the parametrization of $\Grass{2}{4}$ as $(\S^2\times\S^2)/\set{\pm 1}$, which we use in \cref{sec:grassmann-tdesigns} to construct $t$-design curves on $\Grass{2}{4}$ from $t$-design curves and $t$-design points on $\S^2$ .

\subsection{Polynomials and $t$-designs on the Grassmannian}

We denote by $\Pol_t(\Grass{m}{d})$ the collection of classical polynomials of degree at most $t\in\N$ on the Grassmannian $\Grass{m}{d}$:
\begin{equation*}
    \Pol_t(\Grass{m}{d})=\set{f\restricted_{\Grass{m}{d}}\st f\in\R[X]_t},
\end{equation*}
where $\R[X]_t$ is the collection of multivariate polynomials of degree at most $t$ in the $d^2$ variables arranged as a matrix $X\in\R^{d\times d}$.

The eigenspaces of the Laplace--Beltrami operator on $\Grass{m}{d}$ decompose into irreducible subspaces, which are naturally indexed by integer partitions $\pi=(\pi_1,\ldots,\pi_{\ell(\pi)})$ with $\pi_1\geq\cdots\geq \pi_{\ell(\pi)}\geq 0$ and $\ell(\pi)\leq \min\set{m,d-m}$, of \emph{degree} $\DegreePart{\pi}\defined \pi_1+\cdots+\pi_{\ell(\pi)}$. We denote by $H_{\pi}(\Grass{m}{d})$ the irreducible subspace associated with the partition $\pi$, which is equivalent to the irreducible representation $\mathcal{H}^d_{2\pi}$ of $\O(d)$ with type $2\pi$. It is well known (see, for example, \cite[Section 2]{EhlerGraef2019}) that $\Pol_t(\Grass{m}{d})$ can be identified with the following direct sum of these irreducible subspaces:
\begin{equation}\label{eq:Pol-Hpi}
    \Pol_t(\Grass{m}{d})=\bigdirectsum_{\mathclap{\substack{\DegreePart{\pi}\leq t\\ \ell(\pi)\leq \min\set{m,d-m}}}} H_{\pi}(\Grass{m}{d}).
\end{equation}

The notion of a $t$-design extends analogously from the sphere to the Grassmannian. A finite set $X=\set{x_1,\dotsc,x_N}\contained\Grass{m}{d}$ is a \emph{$t$-design} on $\Grass{m}{d}$ if
\begin{equation*}
    \frac{1}{N}\sum_{i=1}^{N}f(x_i)=\int_{\Grass{m}{d}}f\quad\text{for all $f\in\Pol_t(\Grass{m}{d})$}.
\end{equation*}
Designs on Grassmannians were first studied by Bachoc, Coulangeon, and Nebe \cite{tdesignpoints-grassmann-bachoc2002designs} in connection with lattices, and the associated theory of codes and designs in Grassmannian spaces was developed by Bachoc, Bannai, and Coulangeon \cite{tdesignpoints-grassmann-bachoc2004codes} in analogy with the classical theory of Delsarte, Goethals, and Seidel on the sphere. However, the existence of asymptotically optimal $t$-designs on $\Grass{m}{d}$ was settled only indirectly: Etayo, Marzo, and Ortega-Cerdà \cite{tdesignpoints-etayo2018asymptotically} proved the analogue of the Korevaar--Meyers conjecture for $t$-designs on compact algebraic manifolds, which includes the Grassmannians. Their result, specialized to the Grassmannian $\Grass{m}{d}$, is stated in the following theorem.

\begin{theorem}[{\cite[Theorem 2.2]{tdesignpoints-etayo2018asymptotically}}]\label{thm:optimal-points-Grassmannian}
	For $m,d\in\N$ with $m<d$, there is a constant $C_{m,d}$ depending only on $m$ and $d$ such that for each $N\geq C_{m,d}t^{\dim(\Grass{m}{d})}$ there are $t$-designs in $\Grass{m}{d}$ with $N$ points.
\end{theorem}

From now on, we focus on the Grassmannian $\Grass{2}{4}$.

\subsection{Parametrization of $\Grass{2}{4}$ as $(\S^2\times\S^2)/\{\pm 1\}$}

The Grassmannian $\Grass{2}{4}$ is the simplest Grassmannian that is algebraically different from a projective space. It is well-known that $\S^2\times \S^2$ is a double covering of $\Grass{2}{4}$. Specifically, there is an isometric one-to-one map $\ParamGrass\from (\S^2\times\S^2)/\{\pm 1\}\to \Grass{2}{4}$ given by
\begin{equation}\label{eq:paramGrass}
    \ParamGrass(-x,-y)=\ParamGrass(x,y)
    \defined
		\frac{1}{2}\begin{pmatrix}
			1+x\transpose y & -(x\times y)\transpose\\[3pt] -x\times y & xy\transpose+yx\transpose+(1-x\transpose y)\Id_3
		\end{pmatrix},
\end{equation}
where $\Id_3$ is the $3\times 3$ identity matrix and $A\transpose$ denotes the transpose of $A$; see \cite[Section 7]{DickEhlerGraefKrattenthaler2023} for details. Moreover, the pushforward of the uniform probability measure on $\S^2\times \S^2$ under $\ParamGrass$ is the uniform probability measure $\sigma_{\Grass{2}{4}}$ on $\Grass{2}{4}$; see \cite[Lemma F.6]{DickEhlerGraefKrattenthaler2023}. Throughout, unless otherwise stated, all measures are normalized to have total mass $1$.

Recall that the irreducible subspaces of the eigenspaces of the Laplace--Beltrami operator on $\Grass{2}{4}$ are naturally indexed by integer partitions $\pi=(\pi_1,\pi_2)\in\N^2$, with $\pi_1\geq\pi_2\geq 0$, of \emph{degree} $\DegreePart{\pi}\defined \pi_1+\pi_2$; see \cite{JamesConstantine1974} or \cite{tdesignpoints-grassmann-bachoc2002designs} for details. Then, we have the following decomposition of the space $L^2(\Grass{2}{4})$ of square-integrable functions on $\Grass{2}{4}$ into the direct sum of these irreducible subspaces:
\begin{equation*}
    L^2(\Grass{2}{4})=\overline{\bigdirectsum_{\mathclap{\pi_1\geq \pi_2\geq 0}} H_{\pi}(\Grass{2}{4})},
\end{equation*}
where $H_{\pi}(\Grass{2}{4})$ is the irreducible subspace associated with the partition $\pi$; see \cite[(7.1)]{DickEhlerGraefKrattenthaler2023}.

A consequence of the identification of $\Grass{2}{4}$ with $(\S^2\times\S^2)/\{\pm 1\}$ given by \eqref{eq:paramGrass} is that, given an orthonormal basis $\set{\varphi_{\pi,l}}_{1\leq l\leq \dim(H_{\pi}(\Grass{2}{4}))}$ of the irreducible subspace $H_{\pi}(\Grass{2}{4})$, we can express the eigenfunctions $\varphi_{\pi,l}$ as a tensor product of spherical harmonics on $\S^2$. More specifically, consider the functions $Y_{l,m}^{j,k}\from \Grass{2}{4}\to\C$ given by
\begin{equation}\label{eq:tensor-product}
    Y_{l,m}^{j,k}(\ParamGrass(x,y))\defined (Y_{l}^{j}\otimes Y_{m}^{k})(x,y)=Y_{l}^{j}(x)Y_{m}^{k}(y),
\end{equation}
where $Y_{l}^{j}$ denotes a spherical harmonic of degree $l$. These functions are well-defined for $l+m\in2\N$, the latter taking into account the ambiguity in~\eqref{eq:paramGrass}.

\begin{theorem}[{\cite[Theorem 7.2]{DickEhlerGraefKrattenthaler2023}}]\label{thm:span-eigenfunctions}
    Let $H_{\pi}(\Grass{2}{4})$ be the irreducible subspace of the eigenspace of the Laplace--Beltrami operator associated with the partition $\pi=(\pi_1,\pi_2)$. Let $l_{\pi}=\pi_1+\pi_2$ and $m_{\pi}=\pi_1-\pi_2$. Then,
    \begin{equation*}
        H_{\pi}(\Grass{2}{4})=
        \Span\set{Y_{l_{\pi},m_{\pi}}^{j,k},Y_{m_{\pi},l_{\pi}}^{k,j}\st  -l_{\pi}\leq j\leq l_{\pi},
        -m_{\pi}\leq k\leq m_{\pi}}.
    \end{equation*}
\end{theorem}

Next, we define the concept of $t$-design curves on the Grassmannian.

\section{$t$-design curves on the Grassmannian}\label{sec:grassmann-tdesigns}

In contrast to $t$-design points, there are no results about $t$-design curves on the Grassmann manifold so far. Here we focus on the Grassmannian $\Grass{2}{4}$, for which we define $t$-design curves as follows.

\begin{definition}
    For $t\in\N$, we say that a closed, piecewise smooth curve $\gamma\from[0,1]\to\Grass{2}{4}$ is a \emph{$t$-design curve} on $\Grass{2}{4}$ if
    \begin{equation*}
        \frac{1}{\length(\gamma)}\int_{\gamma}f=\int_{\Grass{2}{4}}f\quad \text{for all $f\in\Pol_t(\Grass{2}{4})$}.
    \end{equation*}
\end{definition}

The lower bound on the length of $t$-design curves in \cref{thm:curve-lowerbound} is originally proved for the sphere $\S^d$, but, as already mentioned in \cite{tdesigncurves-EhlerGrochenig2023}, its proof directly carries over to $d$-dimensional compact Riemannian manifolds with the appropriate notion of $t$-designs. Since $\Grass{2}{4}$ is $4$-dimensional, every sequence $(\gamma_t)_{t\in\N}$ of $t$-design curves on $\Grass{2}{4}$ therefore satisfies $\length(\gamma_t)\geq ct^3$ for some constant $c>0$. Thus, we look for sequences of $t$-design curves matching this lower bound.

Our main result is the following theorem, which establishes the existence of asymptotically optimal $t$-design curves on $\Grass{2}{4}$; it can be seen as the analogue of \cref{thm:optimal-curves-S2} for the case of the Grassmann manifold $\Grass{2}{4}$.

\begin{theorem}\label{thm:optimal-curves-G24}
	There exists a sequence $(\Gamma_t)_{t\in\N}$ of $t$-design curves on $\Grass{2}{4}$ with length $\length(\Gamma_t)\asymp t^3$.
\end{theorem}

\begin{proof}
For every $t\in\N$, let $X=\set{x_1,\dotsc,x_n}\contained\S^2$ be a $t$-design with $n\asymp t^2$ points, and let $\gamma$ be a $t$-design curve on $\S^2$ satisfying $\length(\gamma)\asymp t$. We may assume that $\gamma$ contains a pair of antipodal points $p$ and $-p$. Otherwise, let $p\in\gamma$; since $\gamma$ is a $1$-design, there exists a point $a\in\gamma$ with $\innerprod{a}{p}=0$; if $R$ is the rotation by $\pi$ about $a$, the concatenation of $\gamma$ and $R\gamma$ at $a$ is a closed $t$-design curve of length $2\length(\gamma)$ containing $p$ and $-p$. The continuous map $s\mapsto d(p,\gamma(s))$ attains the values $0$ and $\pi$, and hence every value in $[0,\pi]$. As $\S^2$ is two-point homogeneous, that is, $\SO(3)$ acts transitively on pairs of points at a given distance, for all $i,j$ there exists a rotation $O_{ij}\in \SO(3)$ such that $x_i$ and $x_j$ lie on $\gamma^{(i,j)}\defined O_{ij}\gamma$. Consider the curve $\Gamma_t$ on $(\S^2\times\S^2)/\{\pm 1\}$ given by the following set of curve-point pairs on $\S^2\times\S^2$:
\begin{equation*}
    \Gamma_t=\left\{
    \begin{matrix}
        (\gamma^{(1,2)},x_1),  & (x_1,\gamma^{(1,2)}),\\
        (\gamma^{(1,3)},x_2), & (x_2,\gamma^{(1,3)}),\\
        (\gamma^{(2,4)},x_3), & (x_3,\gamma^{(2,4)}),\\
        \vdots & \vdots\\
        (\gamma^{(n-2,n)},x_{n-1}), & (x_{n-1},\gamma^{(n-2,n)}),\\
        (\gamma^{(n-1,n)},x_n), & (x_{n},\gamma^{(n-1,n)})\\
    \end{matrix}
    \right\},
\end{equation*}
that is,
\begin{align*}
    \Gamma_t&=\set[\Big]{(\gamma^{(1,2)},x_1),(x_1,\gamma^{(1,2)}), (\gamma^{(n-1,n)},x_n),(x_n,\gamma^{(n-1,n)})}\\
    &\quad\union\bigunion_{i=2}^{n-1}\set{(\gamma^{(i-1,i+1)},x_i)}\union\bigunion_{i=2}^{n-1} \set{(x_i,\gamma^{(i-1,i+1)})}.
\end{align*}
By construction, $\Gamma_t$ is indeed connected in $(\S^2\times\S^2)/\{\pm 1\}$ and hence induces a connected curve $\mathcal{P}\circ \Gamma_t$ on $\Grass{2}{4}$. Now we must check that it is a $t$-design curve of length $\length(\mathcal{P}\circ \Gamma_t)\asymp t^3$ on $\Grass{2}{4}$. 

According to \eqref{eq:Pol-Hpi}, to prove that $\ParamGrass\composed \Gamma_t$ is a $t$-design curve on $\Grass{2}{4}$, it suffices to check the integration conditions for the eigenfunctions $Y^{j,k}_{l_\pi,m_\pi},Y^{k,j}_{m_\pi,l_\pi}$ in \cref{thm:span-eigenfunctions}. We will perform the computations explicitly for the eigenfunctions of the form $Y^{j,k}_{l_\pi,m_\pi}$; the other case is analogous. The isometric parametrization via $\mathcal{P}$ yields
\begin{align*}
\frac{1}{\ell(\mathcal{P}\circ\Gamma_t)} \int_{\mathcal{P}\circ\Gamma_t} Y^{j,k}_{l_\pi,m_\pi} & = \frac{1}{\length(\Gamma_t)} \int_{\Gamma_t} Y^{j,k}_{l_\pi,m_\pi}\circ\mathcal{P}.
\end{align*}
On the other hand, since $\ParamGrass$ pushes the uniform measure on $\S^2\times\S^2$ forward to $\sigma_{\Grass{2}{4}}$, the tensor structure \eqref{eq:tensor-product} gives
\begin{equation*}
    \int_{\Grass{2}{4}}Y_{l_{\pi},m_{\pi}}^{j,k}\dd\sigma=\int_{\S^2\times\S^2}Y_{l_{\pi},m_{\pi}}^{j,k}\circ \ParamGrass
    =
    \int_{\S^2\times\S^2}Y^j_{l_{\pi}}\otimes Y^k_{m_{\pi}}.
\end{equation*}
Hence, we must verify that 
\begin{equation*}
    \frac{1}{\length(\Gamma_t)}\int_{\Gamma_t}Y_{l_{\pi},m_{\pi}}^{j,k}\circ \ParamGrass
    =
    \int_{\S^2\times\S^2}Y^j_{l_{\pi}}\otimes Y^k_{m_{\pi}}.
\end{equation*}
We now prove that this condition is satisfied. From \cref{eq:tensor-product}, we have
\begin{align*}
    \int_{\Gamma_t}Y_{l_{\pi},m_{\pi}}^{j,k}\circ\ParamGrass &=\int_{\Gamma_t}Y_{l_{\pi}}^{j}\otimes Y_{m_{\pi}}^{k}\\
    &=\int_{\gamma^{(1,2)}\times\set{x_1}}Y_{l_{\pi}}^{j}\otimes Y_{m_{\pi}}^{k}+\int_{\set{x_1}\times \gamma^{(1,2)}}Y_{l_{\pi}}^{j}\otimes Y_{m_{\pi}}^{k}\\
    &\quad+\int_{\gamma^{(n-1,n)}\times\set{x_n}}Y_{l_{\pi}}^{j}\otimes Y_{m_{\pi}}^{k}+\int_{\set{x_n}\times \gamma^{(n-1,n)}}Y_{l_{\pi}}^{j}\otimes Y_{m_{\pi}}^{k}\\
    &\quad+\sum_{i=2}^{n-1}\int_{\gamma^{(i-1,i+1)}\times\set{x_i}}Y_{l_{\pi}}^{j}\otimes Y_{m_{\pi}}^{k}\\
    &\quad+\sum_{i=2}^{n-1}\int_{\set{x_i}\times \gamma^{(i-1,i+1)}}Y_{l_{\pi}}^{j}\otimes Y_{m_{\pi}}^{k}\\
    &=Y_{m_{\pi}}^{k}(x_1)\int_{\gamma^{(1,2)}}Y_{l_{\pi}}^{j}+Y_{l_{\pi}}^{j}(x_1)\int_{\gamma^{(1,2)}}Y_{m_{\pi}}^{k}\\
    &\quad+Y_{m_{\pi}}^{k}(x_n)\int_{\gamma^{(n-1,n)}}Y_{l_{\pi}}^{j}+Y_{l_{\pi}}^{j}(x_n)\int_{\gamma^{(n-1,n)}}Y_{m_{\pi}}^{k}\\
    &\quad+\sum_{i=2}^{n-1}Y_{m_{\pi}}^{k}(x_i)\int_{\gamma^{(i-1,i+1)}}Y_{l_{\pi}}^{j}+\sum_{i=2}^{n-1}Y_{l_{\pi}}^{j}(x_i)\int_{\gamma^{(i-1,i+1)}}Y_{m_{\pi}}^{k}.
\end{align*}
Then, using that $\gamma$ and hence $\gamma^{(i,j)}$ is a $t$-design curve on $\S^2$ with $\length(\gamma)=\length(\gamma^{(i,j)})$, we have
\begin{align*}
    \frac{1}{\length(\gamma)}\int_{\Gamma_t}Y_{l_{\pi},m_{\pi}}^{j,k}\circ \ParamGrass
    &=Y_{m_{\pi}}^{k}(x_1)\int_{\S^2}Y_{l_{\pi}}^{j}+Y_{l_{\pi}}^{j}(x_1)\int_{\S^2}Y_{m_{\pi}}^{k}\\
    &\quad+Y_{m_{\pi}}^{k}(x_n)\int_{\S^2}Y_{l_{\pi}}^{j}+Y_{l_{\pi}}^{j}(x_n)\int_{\S^2}Y_{m_{\pi}}^{k}\\
    &\quad+\sum_{i=2}^{n-1}Y_{m_{\pi}}^{k}(x_i)\int_{\S^2}Y_{l_{\pi}}^{j}+\sum_{i=2}^{n-1}Y_{l_{\pi}}^{j}(x_i)\int_{\S^2}Y_{m_{\pi}}^{k}\\
    &=\sum_{i=1}^{n}Y_{m_{\pi}}^{k}(x_i)\int_{\S^2}Y_{l_{\pi}}^{j}+\sum_{i=1}^{n}Y_{l_{\pi}}^{j}(x_i)\int_{\S^2}Y_{m_{\pi}}^{k}.
\end{align*}   
Since $\set{x_1,\dotsc,x_n}$ is a spherical $t$-design, each of the two sums equals  $n\int_{\S^2}Y_{l_{\pi}}^{j}\int_{\S^2}Y_{m_{\pi}}^{k}$. As $\Gamma_t$ consists of $2n$ pieces of length $\length(\gamma)$, we have $\length(\Gamma_t)=2n\length(\gamma)$, and therefore
\begin{equation*}
    \frac{1}{\length(\Gamma_t)}\int_{\Gamma_t}Y_{l_{\pi},m_{\pi}}^{j,k}\circ\ParamGrass=\int_{\S^2}Y_{l_{\pi}}^{j}\int_{\S^2}Y_{m_{\pi}}^{k}
    =\int_{\S^2\times\S^2}Y_{l_{\pi}}^{j}\otimes Y_{m_{\pi}}^{k},
\end{equation*}
which is the required identity. Finally, by the choice of $X$ and $\gamma$, we have $n\asymp t^2$ and $\length(\gamma)\asymp t$, so that $\length(\ParamGrass\circ\Gamma_t)=\length(\Gamma_t)=2n\length(\gamma)\asymp t^3$, as wanted.
\end{proof}

\section{Example of a $2$-design curve on $\Grass{2}{4}$}\label{sec:grassmann-example}

In this section we present an explicit construction of a $2$-design curve on $\Grass{2}{4}$ using the parametrization of $\Grass{2}{4}$ by $(\S^2\times\S^2)/\{\pm 1\}$ given in~\eqref{eq:paramGrass}.

\begin{figure}[htp]
\centering
\captionsetup{skip=4pt}
\captionsetup[sub]{skip=3pt}
\begin{subfigure}[b]{.49\linewidth}
\centering
\includegraphics[width=3.9cm]{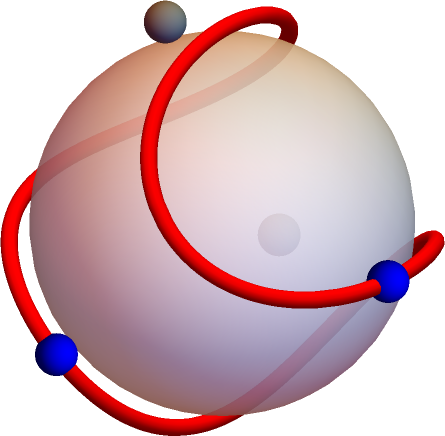}
\caption*{(a1) $\gamma^{(1,2)}$}
\end{subfigure}
\begin{subfigure}[b]{.49\linewidth}
\centering
\includegraphics[width=3.3cm]{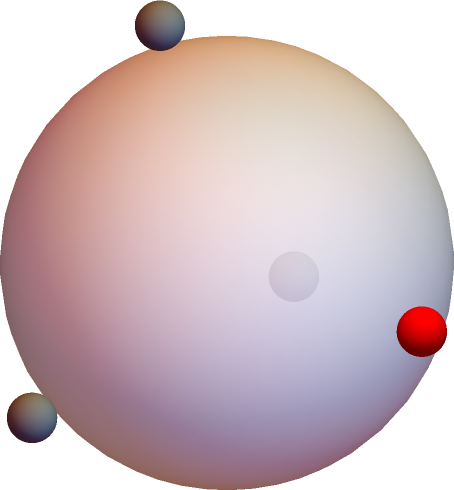}
\caption*{(a2) $x_1$}
\end{subfigure}
\begin{subfigure}[b]{.49\linewidth}
\centering
\includegraphics[width=3.3cm]{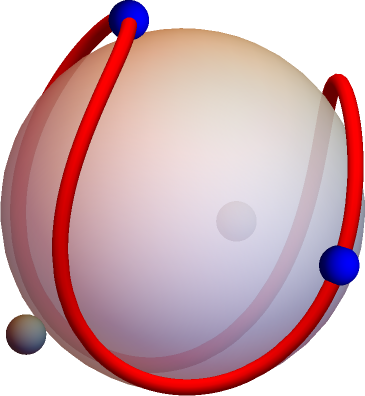}
\caption*{(b1) $\gamma^{(1,3)}$}
\end{subfigure}
\begin{subfigure}[b]{.49\linewidth}
\centering
\includegraphics[width=3.3cm]{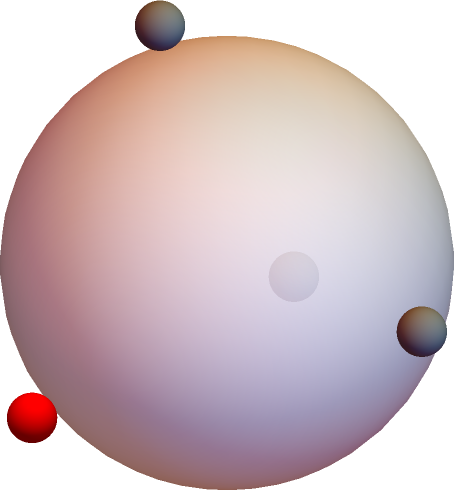}
\caption*{(b2) $x_2$}
\end{subfigure}
\begin{subfigure}[b]{.49\linewidth}
\centering
\includegraphics[width=3.3cm]{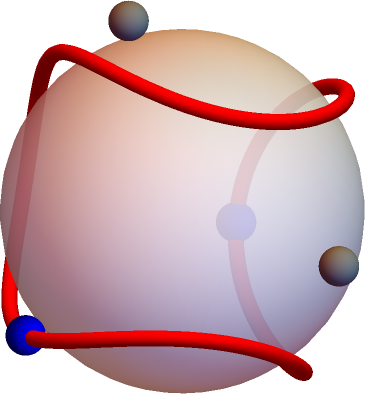}
\caption*{(c1) $\gamma^{(2,4)}$}
\end{subfigure}
\begin{subfigure}[b]{.49\linewidth}
\centering
\includegraphics[width=3.3cm]{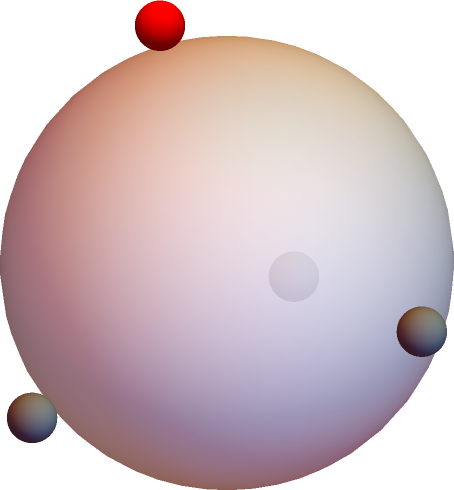}
\caption*{(c2) $x_3$}
\end{subfigure}
\begin{subfigure}[b]{.49\linewidth}
\centering
\includegraphics[width=3.7cm]{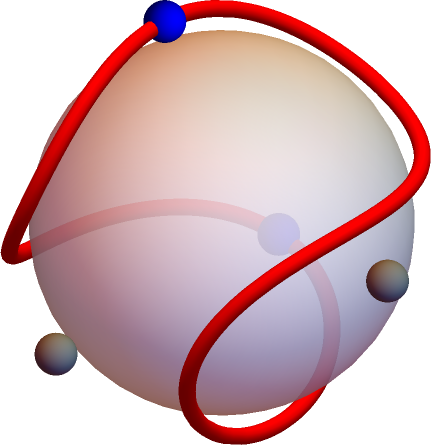}
\caption*{(d1) $\gamma^{(3,4)}$}
\end{subfigure}
\begin{subfigure}[b]{.49\linewidth}
\centering
\includegraphics[width=3.3cm]{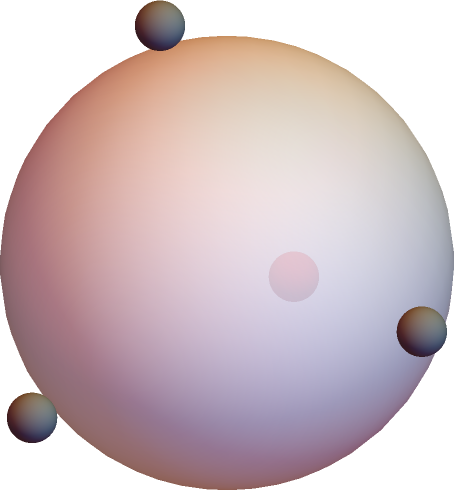}
\caption*{(d2) $x_4$}
\end{subfigure}
    \caption{Visualization of a curve on $\Grass{2}{4}$ using its parametrization as $(\S^2\times\S^2)/\{\pm 1\}$. The image rows form curve-point pairs as in \cref{eq:example-curve}. Each red $2$-design curve $\gamma^{(i,j)}$, $i,j \in \set{1,2,3,4}$, connects the blue $2$-design points $x_i$ and $x_j$ on the left and is associated with the red $2$-design point shown on the right.
    }
    \label{fig:example}
\end{figure}

Consider the vertices of the regular $3$-dimensional simplex, 
\begin{equation*}
    x_1=\begin{pmatrix} \frac{1}{\sqrt{3}}\\\frac{-1}{\sqrt{3}}\\\frac{-1}{\sqrt{3}}\end{pmatrix},\
    x_2=\begin{pmatrix} \frac{-1}{\sqrt{3}}\\\frac{1}{\sqrt{3}}\\\frac{-1}{\sqrt{3}}\end{pmatrix},\
    x_3=\begin{pmatrix} \frac{-1}{\sqrt{3}}\\\frac{-1}{\sqrt{3}}\\\frac{1}{\sqrt{3}}\end{pmatrix},\
    x_4=\begin{pmatrix} \frac{1}{\sqrt{3}}\\\frac{1}{\sqrt{3}}\\\frac{1}{\sqrt{3}}\end{pmatrix},
\end{equation*}
which are $2$-design points on $\S^2$, and the $2$-design curve $\gamma\from [0,1] \to \S^2$ given by
\begin{equation*}
    \gamma(s)=\begin{pmatrix}
        a \cos(2 \pi s)+(1-a)\cos(6\pi s)\\
        a \sin(2 \pi s)-(1-a)\sin(6 \pi s)\\
        2\sqrt{a(1-a)}\sin(4\pi s) 
    \end{pmatrix},
\end{equation*}
where $a \approx 0.7778$ \cite{tdesigncurves-EhlerGrochenig2023}. Since $\gamma(1/4)=(0,1,0)\transpose=-\gamma(3/4)$, we have two antipodal points, and hence there exists a rotation $O_{ij}\in \SO(3)$ such that $x_i$ and $x_j$ lie on $\gamma^{(i,j)}\defined O_{ij}\gamma$.  These points and curves are combined into the curve $\Gamma$ on $\S^2\times\S^2$:
\begin{equation}\label{eq:example-curve}
    \Gamma=\left\{
    \begin{matrix}
        (\gamma^{(1,2)},x_1),  & (x_1,\gamma^{(1,2)}),\\
        (\gamma^{(1,3)},x_2), & (x_2,\gamma^{(1,3)}),\\
        (\gamma^{(2,4)},x_3), & (x_3,\gamma^{(2,4)}),\\
        (\gamma^{(3,4)},x_4), & (x_{4},\gamma^{(3,4)})\\
    \end{matrix}
    \right\}.
\end{equation}
Since $\set{x_1,\dotsc,x_4}$ is a $2$-design and each $\gamma^{(i,j)}$ is a $2$-design curve, the same argument as in the proof of \cref{thm:optimal-curves-G24}, applied with $t=2$ and $n=4$, shows that $\ParamGrass\circ\Gamma$ is indeed a $2$-design curve on $\Grass{2}{4}$.
A visualization of the curves and points on $\S^2$ is provided in \cref{fig:example}.

\section{Harmonic decomposition of the rotation group $\SO(3)$}\label{sec:so3-harmonic}

In this section, we recall the harmonic analysis on the rotation group
\begin{equation*}
    \SO(3)\defined \{R\in \R^{3 \times 3}\st R^{-1}=R\transpose,\ \det(R)=1\}
\end{equation*}
that underlies the definition of polynomials on $\SO(3)$, as well as the identification of $\SO(3)$ with $\S^3/\set{\pm 1}$, which we use in \cref{sec:SO3} to construct $t$-design curves on $\SO(3)$ from curves on $\S^3$. Our presentation follows \cite[Section 2]{graf2012unified}.

\subsection{Polynomials on $\SO(3)$}

For $n\in\N$, we denote by $\PolFixed_n(\S^2)$ the space of spherical harmonics of degree $n$ on the sphere $\S^2$. The irreducible representation $D^n\from \SO(3) \to \mathcal{L}(\PolFixed_n(\S^2))$, where $\mathcal{L}(V)$ denotes the space of linear operators on $V$, of the rotation group $\SO(3)$ on $\PolFixed_n(\S^2)$ is defined as
\begin{equation*}
    D^n(Q)f \defined f(Q\transpose\cdot), \quad f \in \PolFixed_n(\S^2),\ Q \in \SO(3).
\end{equation*}
Choosing a fixed orthonormal basis of spherical harmonics $Y_{n,k}\in \PolFixed_n(\S^2)$, $k=-n,\ldots,n$, yields the matrix representation of the linear operator $D^n(Q)\from\PolFixed_n(\S^2)\to\PolFixed_n(\S^2)$, $Q\in\SO(3)$:
\begin{align*}
    D^n(Q)&\defined (D^n_{k,k'}(Q))_{k,k'=-n,\ldots,n}\in \C^{(2n+1)\times(2n+1)},\\
    D^n_{k,k'}(Q)&\defined\innerprod{D^n(Q)Y_{n,k'}}{Y_{n,k}}_{L^2(\S^2)}=\int_{\S^2}Y_{n,k'}(Q\transpose x)\overline{Y}_{n,k}(x)\dd\mu_{\S^2}(x).
\end{align*}
Note that the indexing of $Y_{n,k}$ differs from that in \cref{sec:grassmann}. The matrix elements $D_{k,k'}^n\from \SO(3) \to \C$, $n\in \N$, $k,k'=-n,\ldots,n$, are smooth functions on $\SO(3)$, called Wigner $D$-functions. By the well-known Peter--Weyl theorem, they form an orthogonal basis of the space $L^2(\SO(3))$.
In analogy with the space $\PolFixed_n(\S^2)$ of spherical harmonics of degree $n$, we define the space of polynomials on $\SO(3)$ of degree exactly $n$ as
\begin{equation*}
    \PolFixed_n(\SO(3))\defined\Span \{D^n_{k,k'}:k,k'=-n,\ldots,n\}, \quad n\in \N.
\end{equation*}

\subsection{Parametrization of $\SO(3)$ as $\S^3/\{\pm 1\}$}

Our strategy for constructing $t$-design curves on $\SO(3)$ in \cref{sec:SO3} is to adapt the asymptotically optimal $t$-design curves on $\S^3$ from \cref{thm:optimal-curves-S3} to the rotation group using that $\S^3$ is a double cover of $\SO(3)$. Recall that points on the sphere $\S^3$ can be written as unit quaternions
\begin{equation*}
    q=\begin{pmatrix} s\\ v \end{pmatrix} \in \S^3, \quad s \in \R, \quad v \in \R^3,
\end{equation*}
with real part $s\in \R$ and vector part $v\in \R^3$. Every unit quaternion $q\in\S^3$ is, up to sign, of the form
\begin{equation}\label{SO3:quaternion_parametrization}
    q(r, \alpha)\defined\begin{pmatrix}
        \cos(\alpha/2)\\
        \sin(\alpha/2)r
    \end{pmatrix}
\end{equation}
for some axis $r\in \S^2$ and some angle $\alpha\in [0,\pi]$. The rotation corresponding to $q(r,\alpha)$ is the rotation $R(r,\alpha)\in \SO(3)$ about the axis $r=(r_1,r_2,r_3)\transpose\in \S^2$, with angle $\alpha \in [0,\pi]$, given by
\begin{equation*}
    R(r,\alpha)=(1-\cos(\alpha)) rr\transpose + \begin{pmatrix}
        \cos(\alpha) & -r_3 \sin(\alpha) & r_2 \sin(\alpha)\\
        r_3 \sin(\alpha) & \cos(\alpha) & - r_1 \sin(\alpha)\\
        -r_2 \sin(\alpha) & r_1 \sin(\alpha) & \cos(\alpha)
    \end{pmatrix}.
\end{equation*}
We can also write this matrix in terms of the coordinates of the quaternion $q\in\S^3$. Writing $q=q(r,\alpha)=(q_1,q_2,q_3,q_4)\transpose$, we have
\begin{equation}\label{eq:Rq}
	R(q)=2\begin{pmatrix}q_2 \\ q_3\\ q_4\end{pmatrix}\begin{pmatrix}q_2 & q_3 &q_4\end{pmatrix}+\begin{pmatrix}
		2q_1^2-1 & -2q_1q_4 & 2q_1q_3\\
		2q_1q_4 & 2q_1^2-1 & -2q_1q_2\\
		-2q_1q_3 & 2q_1q_2 & 2q_1^2-1
	\end{pmatrix}
\end{equation}

We identify the rotation group $\SO(3)$ with the quotient space $\S^3_*\defined\S^3/\{\pm 1\}$ via the isomorphism $q_*\from \SO(3)\to \S^3_*$ given by
\begin{equation}\label{SO3:isomorphism}
    q_*(R(r,\alpha))=\pm q(r,\alpha), \quad r\in \S^2,\ \alpha \in [0,\pi].
\end{equation}
For $x\in\S^3$, we denote the associated point in the quotient space $\S^3_*$ by $x_*\defined\set{x,-x}$. We equip $\S^3_*$ with the quotient metric
\begin{equation*}
    \mathrm{d}_{\S^3_*}(x_{1*},x_{2*}) \defined \min_{\substack{y_1 \in x_{1*} \\ y_2\in x_{2*}}}\mathrm{d}_{\S^3}(y_1,y_2)=\arccos\abs{x_1\transpose x_2}
\end{equation*}
and with the pushforward of the surface measure $\mu_{\S^3}$ under the quotient map $x\mapsto x_*$, that is,
\begin{equation}\label{eq:muS3*}
    \mu_{\S^3_*}(A)\defined\mu_{\S^3}(\set{x\in\S^3\st x_*\in A}),\quad \text{$A\contained\S^3_*$ measurable}.
\end{equation}
On $\SO(3)$, we consider the geodesic distance $\mathrm{d}_{\SO(3)}(O_1,O_2)$ given by the rotation angle of $O_1\transpose O_2$, and the Haar measure $\mu_{\SO(3)}$. We use the normalizations $\mu_{\S^3_*}(\S^3_*)=\mu_{\S^3}(\S^3)= 2\pi^2$ and $\mu_{\SO(3)}(\SO(3))=8 \pi^2$.

\begin{lemma}[{\cite[Lemma 2.1]{graf2012unified}}]\label{lemma:SO3-S3-isomorphism}
    The isomorphism $q_*:\SO(3) \to \S^3_*$ satisfies
    \begin{align}
        \mathrm{d}_{\SO(3)}(O_1,O_2)&=2\mathrm{d}_{\S^3_*}(q_*(O_1),q_*(O_2)), && O_1,O_2 \in \SO(3), \label{dSO3}\\
        \mu_{\SO(3)}(A)&=4\mu_{\S^3_*}(q_*(A)), && \text{$A\contained \SO(3)$ measurable}.\label{muSO3}
    \end{align}
\end{lemma}

By \eqref{dSO3}, the distance between two rotations is twice the distance between the corresponding points of $\S^3_*$; in particular, the image of a curve on $\S^3_*$ under $q_*^{-1}$ has twice the length of the original curve. Similarly, since the map $R\from\S^3\to\SO(3)$ given by $q\mapsto R(q)$ satisfies $R(x)=q_*^{-1}(x_*)$ for $x\in\S^3$, \eqref{muSO3} and \eqref{eq:muS3*} yield, for every integrable function $f$ on $\SO(3)$,
\begin{equation*}
    \int_{\SO(3)}f\dd\mu_{\SO(3)}=4\int_{\S^3}f\circ R\dd\mu_{\S^3},
\end{equation*}
that is,
\begin{equation}\label{eq:integrals-S3-SO3}
    \frac{1}{8\pi^2}\int_{\SO(3)}f\dd\mu_{\SO(3)}=\frac{1}{2\pi^2}\int_{\S^3}f\circ R\dd\mu_{\S^3}.
\end{equation}

In \cref{sec:SO3}, we use these relations to construct asymptotically optimal $t$-design curves on $\SO(3)$ from the $t$-design curves on $\S^3$ provided by \cref{thm:optimal-curves-S3}.

\section{Asymptotically optimal $t$-design curves on $\SO(3)$}\label{sec:SO3}

For $t\in\N$, we denote by $\Pol_t(\SO(3))$ the space of polynomials of degree at most $t$ on the rotation group $\SO(3)$:
\begin{equation}\label{eq:Pol-SO3}
    \Pol_t(\SO(3))\defined\bigdirectsum_{n\leq t}\PolFixed_n(\SO(3)).
\end{equation}
We define $t$-design curves on the rotation group as follows.

\begin{definition}
    For $t\in\N$, we say that a closed, piecewise smooth curve $\Gamma\from[0,1]\to\SO(3)$ is a \emph{$t$-design curve} on $\SO(3)$ if
    \begin{equation*}
        \frac{1}{\length(\Gamma)}\int_{\Gamma}f=\frac{1}{8\pi^2}\int_{\SO(3)}f\,\dd\mu_{\SO(3)}\quad \text{for all $f\in\Pol_t(\SO(3))$}.
    \end{equation*}
\end{definition}
As in the case of the Grassmannian, the lower bound in \cref{thm:curve-lowerbound} carries over to the $3$-dimensional compact Riemannian manifold $\SO(3)$: every sequence $(\Gamma_t)_{t\in\N}$ of $t$-design curves on $\SO(3)$ satisfies
\begin{equation*}
    \length(\Gamma_t)\geq c\,t^2,
\end{equation*}
for some constant $c>0$ \cite{tdesigncurves-EhlerGrochenig2023}.

The previously introduced isomorphism between $\SO(3)$ and $\S^3/\{\pm 1\}$ allows us to generalize \cref{thm:optimal-curves-S3} from $\S^3$ to $\SO(3)$.

\begin{theorem}\label{thm:optimal-curves-SO3}
	There exists a sequence $(\Gamma_t)_{t\in\N}$ of $t$-design curves on $\SO(3)$ with length $\length(\Gamma_t)\asymp t^2$.
\end{theorem}

\begin{proof}
    By \cref{thm:optimal-curves-S3}, let $\gamma_{2t}\from [0,1] \to \S^3$ be a $2t$-design curve with $\length(\gamma_{2t})\asymp t^2$.  Let $g\from \S^3 \to \S^3_*$ be defined as $g(x)=\{x,-x\}$ for $x \in \S^3$.  Note that $q_{*}^{-1}\composed g=R$, where $q_*\from\SO(3)\to\S^3_*$ is the isomorphism given by \eqref{SO3:isomorphism} and $R\from\S^3\to\SO(3)$ is given by \eqref{eq:Rq}. Then, $\Gamma_t\defined R\circ \gamma_{2t}\from [0,1] \to \SO(3)$ is a curve in $\SO(3)$.
    
    Observe that $g$ is a local isometry, as 
    \begin{align*}
        \mathrm{d}_{\S^3}(y_1,y_2)=\mathrm{d}_{\S^3_*}(g(y_{1}),g(y_{2}))
    \end{align*}
    for $\mathrm{d}_{\S^3}(y_1,y_2)<\frac{\pi}{2}$. Furthermore, observe from \cref{dSO3} that
    \begin{align*}
        \mathrm{d}_{\SO(3)}(R(y_1),R(y_2))&=2\mathrm{d}_{\S^3}(y_1,y_2), \quad y_{1},y_{2} \in \S^3,
    \end{align*}
    for $\mathrm{d}_{\S^3}(y_1,y_2)<\frac{\pi}{2}$.
    
    The next step is to prove that $\Gamma_t$ is a $t$-design curve on $\SO(3)$. According to \eqref{eq:Pol-SO3}, it suffices to check the integration property for the Wigner $D$-functions $D_{k,k'}^n$ with $n\leq t$, which form an orthogonal basis for the space $\Pol_t(\SO(3))$. For $x\in\S^2$, let $f_x\from \SO(3)\to\C$ be given by $f_x(Q)\defined Y_{n,k'}(Q\transpose x)$, so that $D_{k,k'}^{n}(Q)=\int_{\S^2}f_x(Q)\overline{Y}_{n,k}(x)\dd\mu_{\S^2}(x)$. By Fubini's theorem,
    \begin{align*}
        \frac{1}{\length(\Gamma_t)}\int_{\Gamma_t}D^n_{k,k'}
        &=\int_{\S^2}\overline{Y}_{n,k}(x)\left(\frac{1}{\length(\Gamma_t)}\int_{\Gamma_t}f_x\right)\dd\mu_{\S^2}(x).
    \end{align*}
    Observe that $\length(\Gamma_t)=2\length(\gamma_{2t})$ and $\int_{\Gamma_t}h=2\int_{\gamma_{2t}}h\circ R$ for every continuous function $h\from\SO(3)\to\C$. Then, we can rewrite the previous expression as
    \begin{align*}
        \frac{1}{\length(\Gamma_t)}\int_{\Gamma_t}D^n_{k,k'}=\int_{\S^2}\overline{Y}_{n,k}(x)\left(\frac{1}{\length(\gamma_{2t})}\int_{\gamma_{2t}}f_x\composed R\right)\dd\mu_{\S^2}(x).
    \end{align*}
    We now use that $\gamma_{2t}$ is a $2t$-design curve: since $f_x$ is a polynomial of degree $n$ in the entries of $Q$ and each entry of $R(q)$ is a quadratic polynomial in the coordinates of $q\in\S^3$, the function $f_x\composed R$ is the restriction to $\S^3$ of a polynomial of degree $2n\leq 2t$. Hence, by the $2t$-design property and \eqref{eq:integrals-S3-SO3},
    \begin{align*}
        \frac{1}{\length(\Gamma_t)}\int_{\Gamma_t}D^n_{k,k'}
        &=\int_{\S^2}\overline{Y}_{n,k}(x)\frac{1}{2\pi^2}\int_{\S^3}(f_x\composed R)(y)\dd\mu_{\S^3}(y)\dd\mu_{\S^2}(x)\\
        &=\int_{\S^2}\overline{Y}_{n,k}(x)\frac{1}{8\pi^2}\int_{\SO(3)}f_x(Q)\dd\mu_{\SO(3)}(Q)\dd\mu_{\S^2}(x)\\
        &=\frac{1}{8\pi^2}\int_{\SO(3)}\int_{\S^2}f_x(Q)\overline{Y}_{n,k}(x)\dd\mu_{\S^2}(x)\dd\mu_{\SO(3)}(Q)\\
        &=\frac{1}{8\pi^2}\int_{\SO(3)}D^n_{k,k'}(Q)\dd\mu_{\SO(3)}(Q).
    \end{align*}
    Thus, we find that $\Gamma_t$ is a $t$-design curve of length $\length(\Gamma_t)=2\length(\gamma_{2t})\asymp t^2$, and by construction there exists a sequence $(\Gamma_t)_{t \in \N}$ of $t$-design curves of asymptotically optimal length.
\end{proof}

\section{Other manifolds}\label{sec:other}

On the torus $\mathbb{T}^d$, Lindblad constructed in \cite{tdesigncurves-lindblad2024asymptotically} $t$-design curves of length of order $t^{d-1}$ with respect to multivariate trigonometric polynomials, and Kunis provided in \cite{kunis2026Td} the corresponding general lower bounds. Following \cite{kunis2026Td}, we denote by
\begin{equation*}
	\Pol_n(\mathbb{T}^d)\defined\Span\set*{x\mapsto e^{2\pi i \innerprod{m}{x}}\st m\in\Z^d,\ \norm{m}_\infty\leq n}
\end{equation*}
the space of multivariate trigonometric polynomials of degree at most $n$ in each variable. For $d,n\in\N$, consider the periodic straight line
\begin{equation*}
	\gamma_{d,n}\from [0,1]\to \mathbb{T}^d,\quad \gamma_{d,n}(y)=\left(y,(2n+1)y,\ldots,(2n+1)^{d-1}y\right)\bmod 1.
\end{equation*}
This is a continuous analogue of a Korobov rank-$1$ lattice rule, and it integrates every $f\in\Pol_{2n}(\mathbb{T}^d)$ exactly, that is,
\begin{equation*}
	\frac{1}{\length(\gamma_{d,n})}\int_{\gamma_{d,n}}f=\int_{\mathbb{T}^d}f.
\end{equation*}
Hence, $\gamma_{d,n}$ is a $t$-design curve on $\mathbb{T}^d$ with $t=2n$, and its length of order $n^{d-1}\asymp t^{d-1}$ is within dimension-dependent constants of the general lower bounds in \cite{kunis2026Td}, so that it is asymptotically optimal.

Regarding the Grassmannian $\Grass{1}{d}$, every element of this space can be viewed as a line through the origin in $\R^d$ and can thus be associated with a pair of antipodal points on the sphere $\S^{d-1}$. Hence, we have
\begin{equation*}
    \Grass{1}{d} \cong \S^{d-1}/\{\pm 1\}.
\end{equation*}
As we did for $\Grass{2}{4}$, we identify $\Grass{1}{d}$ with the set of orthogonal projectors on $\R^d$ of rank $1$ and denote by $\Pol_t(\Grass{1}{d})$ the space of polynomials of degree at most $t$ restricted to $\Grass{1}{d}$. Then, we say that a closed, piecewise smooth curve $\gamma\from[0,1]\to\Grass{1}{d}$ is a \emph{$t$-design curve} on $\Grass{1}{d}$ if
\begin{equation*}
    \frac{1}{\length(\gamma)}\int_{\gamma}f=\int_{\Grass{1}{d}}f\quad \text{for all $f\in\Pol_t(\Grass{1}{d})$},
\end{equation*}
where all measures are normalized to have total mass $1$; in particular, we denote the normalized surface measure on $\S^{d-1}$ by $\sigma_{\S^{d-1}}$ and the uniform measure on $\Grass{1}{d}$ by $\sigma_{\Grass{1}{d}}$.

As mentioned by Lindblad in \cite[Section 3]{tdesigncurves-lindblad2024asymptotically} for the real projective spaces $\RP^2\cong\Grass{1}{3}$ and $\RP^3\cong\Grass{1}{4}$, the previously presented \cref{thm:optimal-curves-S2,thm:optimal-curves-S3} on asymptotically optimal $t$-design curves on $\S^2$ and on $\S^3$ lead to the following result.  Its proof is analogous to that for $\SO(3) \cong \S^3/\{\pm 1\}$ in \cref{sec:SO3}; we include it for completeness.

\begin{proposition}\label{thm:optimal-curves-G1}
    For $d\in\set{3,4}$, there exists a sequence $(\Gamma_t)_{t\in\N}$ of $t$-design curves on $\Grass{1}{d}\cong\RP^{d-1}$ with length $\length(\Gamma_t)\asymp t^{d-2}$.
\end{proposition}

\begin{proof}
Let $d\in\set{3,4}$ and $t\in\N$. By \cref{thm:optimal-curves-S2} for $d=3$ and by \cref{thm:optimal-curves-S3} for $d=4$, there is a $2t$-design curve $\gamma_{2t}\from[0,1]\to\S^{d-1}$ on $\S^{d-1}$ with $\length(\gamma_{2t})\asymp t^{d-2}$. Consider the map
\begin{equation*}
    q\from \S^{d-1}\to\Grass{1}{d},\qquad q(x)\defined xx\transpose,
\end{equation*}
which maps a point $x\in\S^{d-1}$ to the orthogonal projector onto the line spanned by $x$. Since $q(x)=q(y)$ if and only if $y=\pm x$, and every rank-one orthogonal projector is of the form $xx\transpose$ with $x\in\S^{d-1}$, the map $q$ induces a bijection $\S^{d-1}/\{\pm 1\}\to\Grass{1}{d}$; both are models of the real projective space $\RP^{d-1}$. In particular, $q$ is a two-fold covering and scales lengths by the constant factor $\sqrt{2}$. Since $q(O x)=Oq(x)O\transpose$ for all $O\in \O(d)$ and $x\in\S^{d-1}$, the pushforward of $\sigma_{\S^{d-1}}$ under $q$ is an $\O(d)$-invariant probability measure on $\Grass{1}{d}$. By uniqueness of such a measure, it has to be $\sigma_{\Grass{1}{d}}$. Then,
\begin{equation}\label{eq:pushforward-G1}
    \int_{\S^{d-1}}f\composed q\dd\sigma_{\S^{d-1}}=\int_{\Grass{1}{d}}f\dd\sigma_{\Grass{1}{d}}
\end{equation}
for all continuous $f\from\Grass{1}{d}\to\R$. Now, consider the curve
\begin{equation*}
    \Gamma_t\defined q\composed\gamma_{2t}\from[0,1]\to\Grass{1}{d}.
\end{equation*}
Since $\gamma_{2t}$ is closed and piecewise smooth, so is $\Gamma_t$, and its length is $\length(\Gamma_t)=\sqrt{2}\,\length(\gamma_{2t})\asymp t^{d-2}$.

It remains to verify the integration property. Let $f\in\Pol_t(\Grass{1}{d})$. Each entry of $q(x)=xx\transpose$ is a polynomial of degree $2$ in the coordinates of $x$, so the pullback $f\composed q$ is the restriction to $\S^{d-1}$ of a polynomial of degree at most $2t$; moreover, the factor $\sqrt{2}$ cancels, so the normalized line integrals of $f$ along $\Gamma_t$ and of $f\composed q$ along $\gamma_{2t}$ coincide. Since $\gamma_{2t}$ is a $2t$-design curve on $\S^{d-1}$, \cref{eq:pushforward-G1} yields
\begin{equation*}
    \frac{1}{\length(\Gamma_t)}\int_{\Gamma_t}f
    =\frac{1}{\length(\gamma_{2t})}\int_{\gamma_{2t}}f\composed q
    =\int_{\S^{d-1}}f\composed q\dd\sigma_{\S^{d-1}}
    =\int_{\Grass{1}{d}}f\dd\sigma_{\Grass{1}{d}}.
\end{equation*}
Thus, $\Gamma_t$ is a $t$-design curve on $\Grass{1}{d}$ of length $\length(\Gamma_t)\asymp t^{d-2}$.
\end{proof}

Since every subspace of $\R^d$ is uniquely determined by its orthogonal complement, we have
\begin{equation*}
    \Grass{1}{d}\cong\Grass{d-1}{d},
\end{equation*}
where the isometric isomorphism is realized by the complement map
\begin{equation*}
    c\from \Grass{1}{d}\to\Grass{d-1}{d},\qquad c(P)\defined \Id_d-P,
\end{equation*}
which sends a line in $\R^d$ to its orthogonal complement. Since the isometry $c$ preserves lengths and line integrals and maps $\Pol_t(\Grass{d-1}{d})$ to $\Pol_t(\Grass{1}{d})$, applying $c$ to the $t$-design curves on $\Grass{1}{3}$ and $\Grass{1}{4}$ from \cref{thm:optimal-curves-G1} yields the analogous result on $\Grass{2}{3}$ and $\Grass{3}{4}$: there exist sequences $(\Gamma_t)_{t\in\N}$ of $t$-design curves on $\Grass{2}{3}$ with $\length(\Gamma_t)\asymp t$ and on $\Grass{3}{4}$ with $\length(\Gamma_t)\asymp t^2$.

\section*{Acknowledgments}

This work was initiated during a research stay of the third author at the Faculty of Mathematics of the University of Vienna in autumn 2024. He is grateful to the first two authors and to the Faculty of Mathematics for their hospitality.

The authors used the AI systems Claude (Anthropic) and GLM 5.3 (hosted by Mistral) for proofreading and editing the manuscript, as well as for helping with drafting. In a few instances, these systems also suggested corrections to minor inaccuracies or additional justifications for individual steps in the proofs. These suggestions were incorporated after verification by the authors. The mathematical ideas, results, and proofs are due to the authors, who take full responsibility for the content of this paper.

%====================== BIBLIOGRAPHY =====================================

\bibliographystyle{abbrv}
\bibliography{Bibliography.bib}

@Article{tdesigncurves-EhlerGrochenig2023,
  author    = {Ehler, Martin and Gröchenig, Karlheinz},
  journal   = {Forum of Mathematics, Sigma},
  title     = {$t$-Design Curves and Mobile Sampling on the Sphere},
  year      = {2023},
  issn      = {2050-5094},
  pages     = {1--25},
  volume    = {11},
  doi       = {10.1017/fms.2023.106},
  publisher = {Cambridge University Press (CUP)},
}

@article{king2026t2short,
  title={Short spherical $t$-design curves},
  author={King, Emily J. and Mixon, Dustin G.},
  journal={arXiv preprint arXiv:2607.15386},
  year={2026}
}

@Article{ehlergroechenig2026quantitative,
  title={Constraints for stable sampling on the sphere},
  author={Ehler, Martin and Gröchenig, Karlheinz},
  journal={Applied and Computational Harmonic Analysis},
  year={2027},
  volume={86},
  pages = {101922},
  doi = {10.1016/j.acha.2026.101922}
}

@article{tdesigncurves-lindblad2024asymptotically,
  title={Asymptotically optimal $t$-design curves on $\mathbb{S}^3$},
  author={Lindblad, Ayodeji},
  journal={arXiv preprint arXiv:2408.04044v2},
  year={2025}
}

@article{tdesigncurves-lindblad2025generalizations,
  title={Asymptotically short generalizations of $t$-design curves},
  author={Lindblad, Ayodeji},
  journal={arXiv preprint arXiv:2505.03056v3},
  year={2026}
}

@article{kunis2026Td,
  title={Positive quadrature and mobile sampling of multivariate trigonometric polynomials},
  author={Kunis, Stefan},
  journal={arXiv preprint arXiv:2608.11915},
  year={2026}
}

@article{tdesigncurves-ehler2025geodesic,
  title={From point sets to curves: t-designs and {Marcinkiewicz}--{Zygmund} inequalities on the sphere},
  author={Ehler, Martin and Gr{\"o}chenig, Karlheinz and Karner, Clemens},
  journal={Transactions of the American Mathematical Society, Series B},
  volume={13},
  pages={553--578},
  year={2026}
}

@article{tdesignpoints-gariboldi2021optimal,
  title={Optimal asymptotic bounds for designs on manifolds},
  author={Gariboldi, Bianca and Gigante, Giacomo},
  journal={Analysis \& PDE},
  volume={14},
  number={6},
  pages={1701--1724},
  year={2021},
  publisher={Mathematical Sciences Publishers}
}

@article{tdesignpoints-etayo2018asymptotically,
  title={Asymptotically optimal designs on compact algebraic manifolds},
  author={Etayo, Uju{\'e} and Marzo, Jordi and Ortega-Cerdà, Joaquim},
  journal={Monatshefte f{\"u}r Mathematik},
  volume={186},
  pages={235--248},
  year={2018},
  publisher={Springer}
}

@article{tdesignpoints-bondarenko2013optimal,
  author    = {Bondarenko, Andriy and Radchenko, Danylo and Viazovska, Maryna},
  journal   = {Annals of Mathematics},
  title     = {Optimal asymptotic bounds for spherical designs},
  year      = {2013},
  issn      = {0003-486X},
  number    = {2},
  pages     = {443--452},
  volume    = {178},
  doi       = {10.4007/annals.2013.178.2.2},
  publisher = {Annals of Mathematics},
}

@article{tdesignpoints-delsarte1977spherical,
  title={Spherical codes and designs},
  author={Delsarte, Philippe and Goethals, Jean-Marie and Seidel, Johan Jacob},
  journal={Geometriae Dedicata},
  volume={6},
  number={3},
  pages={363--388},
  year={1977},
  publisher={Springer}
}

@article{tdesignpoints-boyvalenkov2001uniqueness,
  title={Uniqueness of the 120-point spherical 11-design in four dimensions},
  author={Boyvalenkov, Peter and Danev, Danyo},
  journal={Archiv der Mathematik},
  volume={77},
  number={4},
  pages={360--368},
  year={2001},
  publisher={Springer}
}

@article{tdesignpoints-goethals1981football,
  title={The football},
  author={Goethals, J. M. and Seidel, J. J.},
  journal={Nieuw Archief voor Wiskunde},
  volume={29},
  number={1},
  pages={50--58},
  year={1981}
}

@article{tdesignpoints-graf2011computation,
  title={On the computation of spherical designs by a new optimization approach based on fast spherical {Fourier} transforms},
  author={Gr{\"a}f, Manuel and Potts, Daniel},
  journal={Numerische Mathematik},
  volume={119},
  number={4},
  pages={699--724},
  year={2011},
  publisher={Springer}
}

@article{tdesignpoints-hardin1996mclaren,
  title={{McLaren’s} improved snub cube and other new spherical designs in three dimensions},
  author={Hardin, Ronald H. and Sloane, Neil J. A.},
  journal={Discrete \& Computational Geometry},
  volume={15},
  pages={429--441},
  year={1996},
  publisher={Springer}
}

@article{tdesignpoints-reznick1995some,
  title={Some constructions of spherical 5-designs},
  author={Reznick, Bruce},
  journal={Linear Algebra and its Applications},
  volume={226},
  pages={163--196},
  year={1995},
  publisher={Elsevier}
}

@article{tdesignpoints-seidel2001definitions,
  title={Definitions for spherical designs},
  author={Seidel, J. J.},
  journal={Journal of Statistical Planning and Inference},
  volume={95},
  number={1-2},
  pages={307--313},
  year={2001},
  publisher={Elsevier}
}

@inproceedings{tdesignpoints-sloane2003spherical,
  title={Spherical designs in four dimensions},
  author={Sloane, Neil J. A. and Hardin, Ronald H. and Cara, Philippe},
  booktitle={Proceedings 2003 {IEEE {Information Theory Workshop (Cat. No. 03EX674)}}},
  pages={253--258},
  year={2003},
  organization={IEEE}
}

@InCollection{tdesignpoints-womersley2018efficient,
  title={Efficient spherical designs with good geometric properties},
  author={Womersley, Robert S.},
  booktitle={Contemporary Computational Mathematics -- A Celebration of the 80th Birthday of Ian Sloan},
  editor={Josef Dick and Frances Y. Kuo and Henryk Woźniakowski},
  pages={1243--1285},
  year={2018},
  publisher={Springer, Cham}
}

@article{tdesignpoints-seymour1984averaging,
  title={Averaging sets: a generalization of mean values and spherical designs},
  author={Seymour, Paul D. and Zaslavsky, Thomas},
  journal={Advances in Mathematics},
  volume={52},
  number={3},
  pages={213--240},
  year={1984},
  publisher={Academic Press}
}

@article{tdesignpoints-korevaar1993spherical,
  title={Spherical {Faraday} cage for the case of equal point charges and {Chebyshev}-type quadrature on the sphere},
  author={Korevaar, Jacob and Meyers, J. L. H.},
  journal={Integral Transforms and Special functions},
  volume={1},
  number={2},
  pages={105--117},
  year={1993},
  publisher={Taylor \& Francis}
}

@article{tdesignpoints-bajnok1992construction,
  title={Construction of spherical t-designs},
  author={Bajnok, Bela},
  journal={Geometriae Dedicata},
  volume={43},
  pages={167--179},
  year={1992},
  publisher={Springer}
}

@article{tdesignpoints-wagner1991averaging,
  title={On averaging sets},
  author={Wagner, Gerold and Volkmann, Bodo},
  journal={Monatshefte f{\"u}r Mathematik},
  volume={111},
  number={1},
  pages={69--78},
  year={1991},
  publisher={Springer}
}

@article{tdesignpoints-ehler2021asymptotically,
  title={Asymptotically optimal cubature formulas on manifolds for prefixed weights},
  author={Ehler, Martin and Etayo, Uju{\'e} and Gariboldi, Bianca and Gigante, Giacomo and Peter, Thomas},
  journal={Journal of Approximation Theory},
  volume={271},
  pages={105632},
  year={2021},
  publisher={Elsevier}
}

@Article{DickEhlerGraefKrattenthaler2023,
  author    = {Josef Dick and Martin Ehler and Manuel Gräf and Christian Krattenthaler},
  journal   = {Constructive Approximation},
  title     = {Spectral Decomposition of Discrepancy Kernels on the {Euclidean} Ball, the Special Orthogonal Group, and the {Grassmannian} Manifold},
  year      = {2023},
  number    = {3},
  pages     = {983--1026},
  volume    = {57},
  doi       = {10.1007/s00365-023-09638-0},
  publisher = {Springer Science and Business Media {LLC}},
}

@article{tdesignpoints-grassmann-bachoc2002designs,
  title={Designs in {Grassmannian} spaces and lattices},
  author={Bachoc, Christine and Coulangeon, Renaud and Nebe, Gabriele},
  journal={Journal of Algebraic Combinatorics},
  volume={16},
  number={1},
  pages={5--19},
  year={2002},
  publisher={Springer}
}

@article{tdesignpoints-grassmann-bachoc2004codes,
  title={Codes and designs in {Grassmannian} spaces},
  author={Bachoc, Christine and Bannai, Eiichi and Coulangeon, Renaud},
  journal={Discrete Mathematics},
  volume={277},
  number={1-3},
  pages={15--28},
  year={2004},
  publisher={Elsevier},
  doi= {10.1016/S0012-365X(03)00151-1},
}

@Article{JamesConstantine1974,
  author       = {Alan T. James and A. G. Constantine},
  title        = {Generalized {Jacobi} Polynomials as Spherical Functions of the {Grassmann} Manifold},
  number       = {1},
  pages        = {174--192},
  volume       = {s3-29},
  year         = {1974},
  doi          = {10.1112/plms/s3-29.1.174},
  journal = {Proceedings of the London Mathematical Society},
  publisher    = {Wiley},
}

@article{applications-breger2024visclust,
  title={{visClust}: A visual clustering algorithm based on orthogonal projections},
  author={Breger, Anna and Karner, Clemens and Ehler, Martin},
  journal={Pattern Recognition},
  volume={148},
  pages={110136},
  year={2024},
  publisher={Elsevier}
}

@article{applications-tSNE,
  author  = {Laurens van der Maaten and Geoffrey Hinton},
  title   = {Visualizing data using {t-SNE}},
  journal = {Journal of Machine Learning Research},
  year    = {2008},
  volume  = {9},
  number  = {86},
  pages   = {2579--2605},
	
}

@Article{EhlerGraefNeumayerSteidl2021,
  author       = {Martin Ehler and Manuel Gräf and Sebastian Neumayer and Gabriele Steidl},
  title        = {Curve Based Approximation of Measures on Manifolds by Discrepancy Minimization},
  number       = {6},
  pages        = {1595--1642},
  volume       = {21},
  year         = {2021},
  doi          = {10.1007/s10208-021-09491-2},
  journal      = {Foundations of Computational Mathematics},
  publisher    = {Springer Science and Business Media {LLC}},
}

@article{graf2012unified,
  title        = {A unified approach to scattered data approximation on {$\mathbb{S}^3$} and {$\operatorname{SO}(3)$}},
  author={Gr{\"a}f, Manuel},
  journal={Advances in Computational Mathematics},
  volume={37},
  number={3},
  pages={379--392},
  year={2012},
  publisher={Springer}
}

@book{applications-so3,
  title={Texture Analysis in Materials Science: Mathematical Methods},
  author={Bunge, H.-J.},
  year={1982},
  publisher={Butterworths, London}
}

@INPROCEEDINGS{sampta-gr24,
  author={Ehler, Martin and Karner, Clemens and López-Gómez, Pedro R.},
  booktitle={2025 International Conference on Sampling Theory and Applications (SampTA)}, 
  title={Asymptotically optimal $t$-design curves on the {Grassmann} manifold}, 
  year={2025},
  volume={},
  number={},
  pages={1-5},
  doi={10.1109/SampTA64769.2025.11133543}}

@article{hielscher2013kernel,
  title={Kernel density estimation on the rotation group and its application to crystallographic texture analysis},
  author={Hielscher, Ralf},
  journal={Journal of Multivariate Analysis},
  volume={119},
  pages={119--143},
  year={2013},
  publisher={Elsevier}
}

@article{hielscher2013funk,
  title={Numerical inversion of the {Funk} transform on the rotation group},
  author={Hielscher, Ralf},
  journal={Inverse Problems},
  volume={29},
  number={12},
  pages={125014},
  year={2013},
  publisher={IOP Publishing}
}

@Article{EhlerGraef2019,
  title={Reproducing kernels for the irreducible components of polynomial spaces on unions of {Grassmannians}},
  author={Ehler, Martin and Gr{\"a}f, Manuel},
  journal={Constructive Approximation},
  volume={49},
  number={1},
  pages={29--58},
  year={2019},
  publisher={Springer}
}
	
\end{document}